# Total Variation Isoperimetric Profiles

Daryl DeFord[*], Hugo Lavenant[†], Zachary Schutzman[‡], and Justin Solomon[§]

**Abstract.** Applications such as political redistricting demand quantitative measures of *geometric compactness* to distinguish between simple and contorted shapes. While the isoperimetric quotient, or ratio of area to perimeter squared, is commonly used in practice, it is sensitive to noisy data and irrelevant geographic features like coastline. These issues addressed in theory by the *isoperimetric profile*, which plots the minimum perimeter needed to inscribe regions of different prescribed areas within the boundary of a shape. Efficient algorithms for computing this profile, however, are not known in practice. Hence, in this paper, we propose a convex Eulerian relaxation of the isoperimetric profile using total variation. We prove theoretical properties of our relaxation, showing that it still satisfies an isoperimetric inequality and yields a convex function of the prescribed area. Furthermore, we provide a discretization of the problem, an optimization technique, and experiments demonstrating the value of our relaxation.

## 1. Introduction and Motivation.

A classic result in modern geometry, the *isoperimetric inequality* states that the least-perimeter shape enclosing a fixed amount of area is a circle. More formally, if $L$ is the length of a simple closed curve in the plane and $A$ is the area it encloses, then $4\pi A \leq L^2$, with equality if and only if the curve is a circle; the observation that a circle minimizes perimeter subject to fixed area is a direct byproduct of this expression. Inspired by this inequality, the *isoperimetric quotient* $4\pi A/L^2$ is a commonly-used proxy for measuring the compactness of a shape; here, *compactness* refers to geometric regularity, a piece of terminology that originated in the political geography literature (see §1.2). This quantity is scale-free, unitless, and bounded between zero and one. It is intended to capture how efficiently a shape "uses" its perimeter to enclose its area and is maximized for a circle.

This quotient, however, is *unstable* in the sense that a small perturbation of the shape's boundary can greatly increase its perimeter without significantly affecting its area; Figure 1 shows one example of this instability. Moreover, as a measure of compactness, the isoperimetric quotient conflates multiple scales. At a coarse scale, a perturbed circle still appears fairly compact, while at a fine scale, the noisy boundary differentiates it from a proper circle; this distinction cannot be captured by a single number.

To address these issues, in this article we investigate a *multiscale* means of measuring compactness that explicitly assesses aspects of compactness at different scales. Our measure is easily and reliably computed, even on very distorted shapes. Although our work on this problem was inspired by concrete challenges related to the evaluation of compactness of voting districts (see §1.2), we find that the isoperimetric profile and its convex relaxation have

---

[*]Massachusetts Institute of Technology
[†]Université Paris-Sud
[‡]University of Pennsylvania
[§]Massachusetts Institute of Technology (corresponding author, jsolomon@mit.edu; authors are alphabetized)

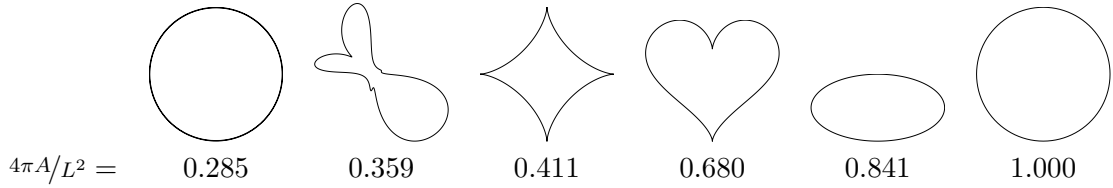

$4\pi A/L^2 =$    0.285      0.359      0.411      0.680      0.841      1.000

Figure 1: A variety of shapes marked with isoperimetric ratios $4\pi A/L^2$ marked; the least compact shape is a small perturbation of a circle.

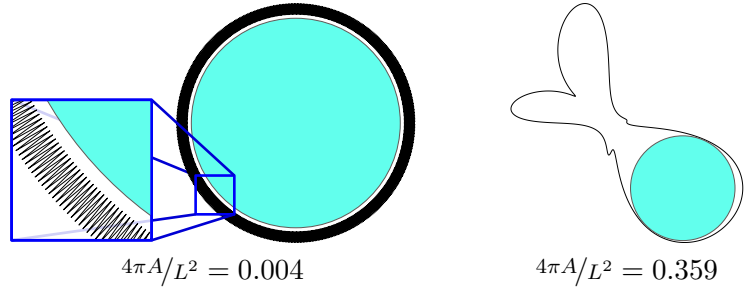

$4\pi A/L^2 = 0.004$          $4\pi A/L^2 = 0.359$

Figure 2: The perturbed circle (left) not compact as measured by the isoperimetric ratio, but it admits a compact inscribed circle with slightly less area; the non-compact shape on the right does not contain a compact inscribed shape that uses most of the interior area.

intrinsic mathematical interest that inspire additional challenging problems in geometry and convex analysis.

### 1.1. Mathematical Overview.

A modern construction in geometry offers a potential resolution to the instability of the isoperimetric ratio. The *isoperimetric profile* of a shape replaces the isoperimetric quotient with a plot of compactness values at different length or area scales. For each area $t$ between 0 and the area $A$ of a shape ($0 \le t \le A$), the isoperimetric profile measures the minimum amount of perimeter *inscribed within* the original shape needed to enclose area $t$. As shown in Figure 2, by considering $t < A$ we can use this to separate perturbative boundary effects from larger-scale barriers to compactness.

Compared to most of the mathematical literature (see e.g. [40]), we consider the *full* perimeter of the inscribed shape rather than the relative perimeter; the latter does not account for perimeter shared with the boundary of the input shape. Considering only relative perimeters—or, equivalently, the isoperimetric problem in domains without boundaries—leads to an elegant mathematical structure with a strong interplay between the geometry of the domain and the shape of the profile. The problem considered in our work, which takes the full perimeter into account, is less studied theoretically (see [45] for some relevant discussion, as well as the forthcoming article [44] which was informed by a preprint of the present article), likely because the shape of the profile has a less elegant qualitative properties. This case, however, is relevant for computational problems in need of stable, multiscale compactness measures.

Although the isoperimetric profile is a promising theoretical construction, to-date no algorithms are known to compute this plot in practice. While special properties in two dimensions suggest that a computational geometry technique *may* be possible to formulate for polygons in the plane, this problem remains open; even less is known about computing the profile in the higher-dimensional case or on Riemannian manifolds. One possible computational approach might use a phase field approximation *à la* Modica–Mortola [10], but the resulting problem would require the minimization of a nonconvex functional.

Given the potential applications of the isoperimetric profile and the challenges of computing it in practice, this paper proposes a convex relaxation built from total variation. Our basic approach is to write the optimization problem underlying the isoperimetric profile in Eulerian language and then to relax an integer variable to be real-valued. The resulting problem is straightforward to optimize after using standard discretization techniques from mathematical imaging. We show that a key theoretical property of the isoperimetric profile is preserved in our convex relaxation, namely that the lower envelope of our set of relaxed profiles is provided by a circle—the most compact shape; we also provide theoretical results giving qualitative intuition for the behavior of our profile and its relationship to the original nonconvex problem.

### 1.2. Compactness and Political Redistricting.
A key application of the isoperimetric quotient—known as the *Polsby–Popper* measure in political science [36]—is in measuring the compactness of legislative voting districts. Here, we use the word *compact* as it used in the political geography literature, to refer to plane regions that appear sufficiently regular, as there is no legally agreed-upon definition of this concept; see [29] for further analysis and discussion.

In representative systems like the U.S. House of Representatives, voters are clustered geographically into districts with each district electing a single representative to Congress. Manipulating district shapes to engineer a particular outcome for a vote is a practice known as *gerrymandering*, which undermines democratic principles and has been used to deny underrepresented minorities the opportunity to elect a representative of their choosing. While a specific measure of geometric quality is not always mandated by law, scores like Polsby–Popper are used as quantitative proxies for the reasonableness of a district or districting plan; contorted, nonconvex district shapes may signal that a district was designed with motives other than those required by law.

While compactness scores like the isoperimetric quotient are widely used in arguments for or against districting plans, they have significant mathematical drawbacks that undermine their interpretability and reliability [21]. Most importantly, these scores are *unstable.* As shown in Figure 1, the Polsby–Popper score is unstable under boundary perturbations. A visually insignificant adjustment to the boundary of a circle, which has isoperimetric quotient equal to one, can make the isoperimetric quotient arbitrarily close to zero; a minor perturbation changes the shape from being the 'most' compact to the 'least' compact. This is by no means a degenerate edge-case, but rather extremely common in geographic information systems (GIS) data due to fractal-like structures resulting from the degree of precision with which the geography is measured; maps of different *resolution* of the same district can lead to Polsby–Popper scores that differ by a factor of $2-3$ [3]. This particular sensitivity can greatly distort the evaluation of the compactness of a district, since natural choices of boundaries for districts, such as coastlines, geological features, or municipality boundaries may be noisy.

This instability makes it difficult for the isoperimetric quotient to distinguish between districts whose boundaries are contorted to gerrymander rather than those that simply contain a coastline or follow a complicated municipal boundary.

The isoperimetric profile, which replaces a single isoperimetric quotient with a whole plot of compactness scores and whose convex relaxation can be computed reliably, can be used to discriminate between compact and noncompact districts without these drawbacks. To demonstrate the interest of this new measurement, in this paper we accompany synthetic experiments with illustrations of the (relaxed) isoperimetric profile of actual U.S. voting districts in the state of North Carolina.

We mention a related notion in the political science literature: A compactness score measuring the proportion of the area of the region filled by its largest inscribed circle is known as the *Ehrenburg test*, presented in a critique of various compactness measures developed in the nineteenth century to describe the geographic and human features of landmasses [22]. This measure captures information about an intermediate scale that depends on the shape of the region being considered but again summarizes compactness with a single number.

## 2. Background and Preliminaries.
We begin our discussion with some basic notation and background information setting the stage for our mathematical construction in §3.

### 2.1. Isoperimetry.
Although the basic isoperimetric problem dates back centuries, the development of isoperimetric inequalities remains an active area of research in mathematics; for a general survey, we refer the reader to [40] for a broad discussion of classical results and open problems. Here, we restrict to basic Euclidean constructions relevant to our computational task.

Take $\Omega \subseteq \mathbb{R}^n$ to be a compact region whose boundary $\partial \Omega \subseteq \Omega$ is an $(n-1)$-dimensional hypersurface in $\mathbb{R}^n$, and let $B_1$ be the unit ball in that space. The general isoperimetric inequality states

$$(2.1) \qquad n \cdot \mathrm{vol}(\Omega)^{(n-1)/n} \cdot \mathrm{vol}(B_1)^{1/n} \leq \mathrm{area}(\partial \Omega).$$

This expression encodes the fact that the unit ball $B_1$ minimizes $(n-1)$-dimensional boundary surface area over all regions $\Omega$ in $\mathbb{R}^n$ with volume 1. As we have already seen, however, it is straightforward to perturb the boundary of $\Omega$ slightly to increase the right-hand side of this expression an arbitrary amount with minor impact on the left-hand side; this is a potential source of instability in practice.

More generally, take $\Omega \subseteq \mathbb{R}^n$ to be the same region, and take $t \in [0, \mathrm{vol}(\Omega)]$. For each $t$ we could ask for the smallest surface area needed to enclose volume $t$ completely within $\Omega$:

$$(2.2) \qquad I_\Omega(t) := \min\{\mathrm{area}(\partial \Sigma) : \Sigma \subseteq \Omega \text{ and } \mathrm{vol}(\Sigma) = t\}.$$

Here, we define $\mathrm{area}(\partial \Sigma)$ to be the area of $\partial \Sigma$ as a submanifold of $\mathbb{R}^n$, that is, including the area of the intersection $\partial \Sigma \cap \partial \Omega$; we put no topological restrictions—in particular, connectedness—on $\Sigma$. For our analysis in later sections to be relevant, the minimal regularity assumption needed on $\Omega$ is that the TV profile $I_\Omega^{\mathrm{TV}}(t)$ (to be defined in (3.2)) is finite whenever $t < \mathrm{vol}(\Omega)$. This condition can be checked with minimal assumptions on $\partial \Omega$ (including fractality); see Appendix A for discussion. A direct consequence of the isoperimetric inequality (2.1) is that for any fixed $t$ and prescribed volume $V$, minimizers of $I_\Omega(t)$ over all possible regions $\Omega \subseteq \mathbb{R}^n$ with $\mathrm{vol}(\Omega) = V$ are those containing a ball with volume $t$.

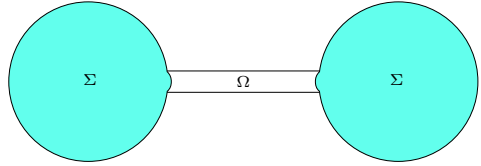

Figure 3: The isoperimetric profile of a nonconvex shape may require a disconnected domain $\Sigma$ inscribed in $\Omega$; for example, the most efficient use of area inscribed in the barbell shape above might be disconnected since the rectangular portion in the middle has little area but considerable perimeter. Note this is a schematic rather than computed algorithmically.

In any event, under fairly weak assumptions, a minimizer of problem (2.2) is known to exist. Away from $\partial\Omega$, the boundary $\partial\Sigma$ has constant mean curvature and is nonsingular for $n < 8$ [24]. Because we include the area of $\partial\Sigma \cap \partial\Omega$ in our variational problem, $\partial\Sigma$ and $\partial\Omega$ meet tangentially. The free boundaries of $\Sigma$ are minimal surfaces (see e.g. [40, Theorem 1]).

When $n = 2$, the free boundaries of $\Sigma$ are circular arcs with the same signed radius [5, 6]. This property *may* be helpful in designing a computational geometry-style algorithm for computing $I_\Omega(t)$ when $\Omega$ is a polygon in $\mathbb{R}^2$, although to our knowledge no such algorithm has been proposed; it remains unknown in which cases computing $I_\Omega(t)$ can be done in polynomial time. As shown in Figure 3, the optimal $\Sigma$ might be disconnected, even in 2D. Additional open problems in the $n = 2$ case are posed in [20]. To make matters more challenging, this helpful structure does not appear to extend for $n \geq 3$, underscoring the potential value of our construction in §3 even if $I_\Omega$ is efficiently computable in some special cases.

**Open Problem 2.1.** *Identify a polynomial-time algorithm or NP-hardness result for computing isoperimetric profiles. The simplest open problem is computing the isoperimetric profile of a polygon in the plane $\mathbb{R}^2$; if evaluating the profile of such a polygon is polynomial-time solvable, higher-dimensional analogs are of interest as well.*

## 2.2. Total Variation.

The definition of the isoperimetric profile (2.2) is Lagrangian: The variable $\Sigma$ explicitly parameterizes the boundary of an unknown shape. Our approach in §3 will be to switch an Eulerian formulation, replacing the boundary shape optimization variable with the indicator function of the unknown $\Sigma$. To do so, we need an Eulerian way to compute the area of $\partial\Sigma$, which we can achieve using the *total variation* of its indicator.

For a function $f \in L^1(\mathbb{R}^n)$, the total variation (TV) of $f$ is defined as

$$(2.3) \qquad \mathrm{TV}[f] := \sup \left\{ \int_{\mathbb{R}^n} [f(x)\nabla \cdot \phi(x)]\, dx : \phi \in C_c^1(\mathbb{R}^n \to \mathbb{R}^n) \text{ and } \|\phi\|_\infty \leq 1 \right\}.$$

Here, $C_c^1(\mathbb{R}^n \to \mathbb{R}^n)$ denotes the space of compactly-supported continuously differentiable functions from $\mathbb{R}^n$ into $\mathbb{R}^n$. For differentiable functions $f \in C^1(\mathbb{R}^n)$ with compact support, the divergence theorem implies an alternative formula for total variation introduced to math-ematical image processing in [42]:

$$(2.4) \qquad \mathrm{TV}[f] = \int_{\mathbb{R}^n} \|\nabla f\|_2 \, dx.$$

An alternative way of computing total variation, which reveals its link with the isoperimetric problem, is the co-area formula [23]. For a non-negative $f$, this formula states

$$(2.5) \qquad \mathrm{TV}[f] = \int_0^{+\infty} \mathrm{area}(\partial\{f \geq s\}) ds.$$

For a region $\Sigma \subseteq \mathbb{R}^n$, denote its *indicator function* $\mathbb{1}_\Sigma$ via

$$(2.6) \qquad \mathbb{1}_\Sigma(x) := \left\{ \begin{array}{ll} 1 & \text{if } x \in \Sigma \\ 0 & \text{otherwise.} \end{array} \right.$$

Then, a consequence of the co-area formula (2.5) is that

$$(2.7) \qquad \mathrm{area}(\partial\Sigma) = \mathrm{TV}[\mathbb{1}_\Sigma].$$

A detailed account of geometric information measured by total variation is provided in [15].

**3. Total Variation Isoperimetric Profile.** Now that we have established basic notation and the functionals we will use in our construction, we are ready to present our convex relaxation of the isoperimetric profile for a shape embedded in $\mathbb{R}^n$.

**3.1. Definition.** Inspired by the total variation formula (2.7), we can formally rewrite the optimization problem (2.2) for the isoperimetric profile in an Eulerian fashion:

$$(3.1) \qquad I_\Omega(t) = \left\{ \begin{array}{ll} \inf_{f \in L^1(\mathbb{R}^n)} & \mathrm{TV}[f] \\ \text{subject to} & \int_{\mathbb{R}^n} f(x) \, dx = t \\ & 0 \leq f \leq \mathbb{1}_\Omega \\ & f(x) \in \{0, 1\} \; \forall x \in \mathbb{R}^n. \end{array} \right.$$

Here, $f$ is the indicator of the unknown shape $\Sigma$. The three constraints ensure that (1) the area enclosed in $\Sigma$ is equal to $t$, (2) $\Sigma$ is inscribed within $\Omega$, and (3) that $f$ is properly an indicator function. For convenience, we have written this problem using integrals over $\mathbb{R}^n$ but in reality the second constraint ensures that $f$ is zero outside $\Omega$.

Since total variation is convex, the form (3.1) hides all the nonconvexity of the problem in the third constraint. Hence, we propose using an alternative measure of compactness that drops the nonconvex constraint, which we call the *total variation (TV) profile* of $\Omega$:

$$(3.2) \qquad \boxed{I_\Omega^{\mathrm{TV}}(t) := \left\{ \begin{array}{ll} \min_{f \in L^1(\mathbb{R}^n)} & \mathrm{TV}[f] \\ \text{subject to} & \int_{\mathbb{R}^n} f(x) \, dx = t \\ & 0 \leq f \leq \mathbb{1}_\Omega. \end{array} \right.}$$

The existence of a minimizer is an immediate consequence of Rellich's compactness theorem [1, Theorem 3.23] (also see [14, Theorem 2]).

Since (3.2) is obtained by dropping a constraint from (3.1), we immediately find:

**Proposition 3.1.** *For all $(\Omega, t)$, we have $I_\Omega^{\mathrm{TV}}(t) \leq I_\Omega(t)$.*

Unfortunately, this inequality is not tight, which we demonstrate in the following example:

**Example 3.2 (Circle).** *Suppose $\Omega \subset \mathbb{R}^2$ is a circle of radius $R$, and take $t = \pi r^2$ for $r \in (0, R)$. In this case, by the isoperimetric inequality we know $I_\Omega(t) = 2\pi r$. But suppose we take $f(x) \equiv r^2/R^2$. Notice $f(x)$ satisfies the contraints in (3.2), but by the co-area formula*

$$I_\Omega^{\mathrm{TV}}(t) \leq \mathrm{TV}[f] = 2\pi R \cdot \frac{r^2}{R^2} = 2\pi r \cdot \frac{r}{R} < I_\Omega(t).$$

*Hence, our relaxation is not tight.*

As we will prove in Proposition 3.5, the TV profile still admits an elegant characterization explaining why the relaxation is not tight: it is the lower convex envelope of the isoperimetric profile.

The following monotonicity property can be derived directly from the definition of the total variation profile:

**Proposition 3.3.** *Suppose $\Omega_1 \subseteq \Omega_2 \subset \mathbb{R}^d$ are two compact domains. Then for all $t \leq \mathrm{vol}(\Omega_1)$ it holds $I_{\Omega_2}^{\mathrm{TV}}(t) \leq I_{\Omega_1}^{\mathrm{TV}}(t)$.*

In particular, if a domain $\Omega$ can be squeezed between two other domains: $\Omega_1 \subseteq \Omega \subseteq \Omega_2$ then the total variation profile of $\Omega$ is squeezed between those of $\Omega_1$ and $\Omega_2$, at least for $t \leq \mathrm{vol}(\Omega_1)$. For instance, the total variation profile of the perturbed circle (left of Figure 2) must be very close, at least for $t$ away from $\mathrm{vol}(\Omega)$, to a straight line, the total variation profile of a circle. More generally, if we know the domain $\Omega$ with some uncertainty or only at a given resolution, which is roughly equivalent to saying that we can squeeze it between two other domains with a controlled difference in area, then we have guarantees for the value of the TV profile as soon as $t$ is not too close to $\mathrm{vol}(\Omega)$.

**3.2. Duality and convexity.** We next derive a dual for the TV profile problem (3.2) to provide additional insight into the structure of the problem and in particular to show that $I_\Omega^{\mathrm{TV}}(t)$ is convex in $t$. We begin by writing a minimax formulation of the problem by substituting (2.3):

$$(3.3) \qquad I_\Omega^{\mathrm{TV}}(t) = \min_{f \in L^1(\mathbb{R}^n)} \left\{ \begin{array}{ll} \sup_{\phi, \lambda, \psi, \xi} & \int_{\mathbb{R}^n} [f(x) \nabla \cdot \phi(x)] \, dx \\ & + \lambda \left( t - \int_{\mathbb{R}^n} f(x) \, dx \right) \\ & - \int_{\mathbb{R}^n} [\psi(x) f(x) + \xi(x) (\mathbb{1}_\Omega(x) - f(x))] \, dx \\ \text{subject to} & \|\phi\|_\infty \leq 1, \psi \geq 0, \xi \geq 0 \end{array} \right.$$

Here, the dual variables are $\phi \in C_c^1(\mathbb{R}^n \to \mathbb{R}^n)$, $\lambda \in \mathbb{R}$, and $\psi, \xi \in C_c^1(\mathbb{R}^n \to \mathbb{R})$.

A standard argument using the Fenchel–Rockafellar duality theorem [12, Theorem 1.12] justifies swapping the inner and outer problems; for completeness, a formal argument is in Appendix B. Removing extraneous terms, we are left with the following inner variational problem for $f$, with the dual variables fixed:

$$\inf_{f \in L^1(\mathbb{R}^n)} \int_{\mathbb{R}^n} f(x) [\nabla \cdot \phi(x) - \lambda - \psi(x) + \xi(x)] \, dx.$$

This problem is unbounded unless $\nabla \cdot \phi(x) + \xi(x) - \psi(x) = \lambda$ a.e. $x \in \mathbb{R}^n$. Hence, after swapping the inner and outer problems in (3.3) and substituting this relationship we find

$$I_\Omega^{\mathrm{TV}}(t) = \begin{cases} \sup_{\phi,\lambda,\psi,\xi} & \lambda t - \int_\Omega \xi(x) \, dx \\ \text{subject to} & \|\phi\|_\infty \leq 1, \psi \geq 0, \xi \geq 0 \\ & \nabla \cdot \phi(x) + \xi(x) - \psi(x) = \lambda \text{ a.e. } x \in \mathbb{R}^n \end{cases}$$

Define $\eta(x) := \xi(x) - \psi(x)$. From the objective of the problem above, we can see $\xi(x) = \max(\eta(x), 0)$ and $\psi(x) = \max(-\eta(x), 0)$. Substituting and simplifying leads to the dual

$$(3.4) \qquad \boxed{I_\Omega^{\mathrm{TV}}(t) = \begin{cases} \sup_{\phi \in C_c^1(\mathbb{R}^n \to \mathbb{R}^n), \lambda \in \mathbb{R}} & \lambda t - \int_\Omega \max(\lambda - \nabla \cdot \phi(x), 0) \, dx \\ \text{subject to} & \|\phi\|_\infty \leq 1 \end{cases}}$$

While this dual formula provides an interesting—if abstract—reinterpretation of $I_\Omega^{\mathrm{TV}}(t)$ in its own right, it also allows us to derive another property of the function itself:

**Proposition 3.4.** *$I_\Omega^{\mathrm{TV}}(t)$ is a convex function of $t$.*

*Proof.* Define
$$h(\lambda) := \inf_{\|\phi\|_\infty \leq 1} \int_\Omega \max(\lambda - \nabla \cdot \phi(x), 0) \, dx.$$

From (3.4), we see that $I_\Omega^{\mathrm{TV}}(t) = h^*(t)$, the convex conjugate (Legendre–Fenchel transform) of $h$, which is necessarily convex [7]. ∎

We can say even more: The TV profile $I_\Omega^{\mathrm{TV}}$ is the lower convex envelope of the isoperimetric profile $I_\Omega$, i.e., the largest convex function upper-bounded by the isoperimetric profile.

**Proposition 3.5.** *The function $I_\Omega^{\mathrm{TV}}$ is the lower convex envelope of $I_\Omega$.*

*Proof.* It is known that the lower convex envelope of a function is equal to the Legendre transform of its Legendre transform [39, Section 12]. Since Proposition 3.4 shows that $I_\Omega^{\mathrm{TV}}$ is convex—and thus equals its convex envelope—it suffices to prove that the Legendre transform of $I_\Omega$ coincides with that of $I_\Omega^{\mathrm{TV}}$.

By definition,

$$(3.5) \qquad (I_\Omega)^\star(\lambda) = \sup_t [\lambda t - I_\Omega(t)] = \begin{cases} \sup_f & \lambda \int_\Omega f(x) dx - \mathrm{TV}(f) \\ \text{subject to} & 0 \leq f \leq \mathbb{1}_\Omega, \\ & f(x) \in \{0, 1\} \ \forall x \in \mathbb{R}^n. \end{cases}$$

On the other hand, the Legendre transform of $I_\Omega^{\mathrm{TV}}$ is given by

$$(3.6) \qquad (I_\Omega^{\mathrm{TV}})^\star(\lambda) = \sup_t \lambda t - I_\Omega(t) = \begin{cases} \sup_f & \lambda \int_\Omega f(x) dx - \mathrm{TV}(f) \\ \text{subject to} & 0 \leq f \leq \mathbb{1}_\Omega. \end{cases}$$

Since there are more admissible competitors in (3.6) than (3.5), $(I_\Omega^{\mathrm{TV}})^\star(\lambda) \geq (I_\Omega)^\star(\lambda)$. On the other hand, using the co-area formula (2.5), for any admissible competitor $f$ in (3.6),

$$\lambda \int_\Omega f(x) dx - \mathrm{TV}(f) = \int_0^1 \underbrace{(\lambda \mathrm{vol}(\{f \geq s\}) - \mathrm{area}(\partial\{f \geq s\}))}_{\leq (I_\Omega)^\star(\lambda) \text{ by definition}} \, ds \leq (I_\Omega)^\star(\lambda).$$

Taking the supremum in $f$ leads to $(I_\Omega^{\mathrm{TV}})^\star(\lambda) \le (I_\Omega)^\star(\lambda)$, which concludes the proof. ∎

The proof of the equality of (3.5) and (3.6) is just a slight adaptation of [14, Proposition 2.1], in a simpler case since we want only the equality of the values of the minimization problem rather than a correspondence between the minimizers.

### 3.3. Behavior for small $t > 0$.

Continuing in our effort to describe the shape of $I_\Omega^{\mathrm{TV}}(t)$ as precisely as possible, we can accompany the convexity result in Proposition 3.4 with a description of its behavior when $t$ is close to zero. In particular, we will show that $I_\Omega^{\mathrm{TV}}$ is *linear* with positive slope for small $t$, a property we can verify in our experiments. This is not surprising given Proposition 3.5, since the isoperimetric profile of $\Omega$ coincides with the one of $\mathbb{R}^n$ for small $t$, the latter being concave. Nevertheless, we can say more by describing exactly the slope of $I_\Omega^{\mathrm{TV}}$ for small volumes.

We use the following auxiliary problem, known as the Cheeger problem [35]. The Cheeger constant of $\Omega$, denoted by $h_1(\Omega)$, is defined as

$$h_1(\Omega) := \inf_{\tilde{\Sigma} \subseteq \Omega} \frac{\mathrm{area}(\partial\tilde{\Sigma})}{\mathrm{vol}(\tilde{\Sigma})},$$

and a subset $\Sigma \subseteq \Omega$ such that $h_1(\Omega) = \frac{\mathrm{area}(\partial\Sigma)}{\mathrm{vol}(\Sigma)}$ is known as a Cheeger set of $\Omega$. Such a set exists as soon as $\Omega$ has a Lipschitz boundary [35, Proposition 3.1] and is unique if $\Omega$ is convex [35, Proposition 5.2]. $\Sigma$ can be interpreted as the largest "smooth" subset of $\Omega$. An explicit description of Cheeger sets can be found for convex domains [30], as well as domains without necks [33]; we emphasize, however, that domains appearing in redistricting applications are not likely to satisfy the assumptions of these articles.

Furthermore, provided we know a Cheeger set of $\Omega$, we can describe the behavior of $I_\Omega^{\mathrm{TV}}(t)$ for small $t$:

**Proposition 3.6.** *Let $\Omega$ be compact, let $h_1(\Omega)$ be the Cheeger constant of $\Omega$, and let $\Sigma$ be a Cheeger set of $\Omega$. Then for any $t \le \mathrm{vol}(\Sigma)$, we have $I_\Omega^{\mathrm{TV}}(t) = h_1(\Omega)t$, and a solution $f$ in (3.2) is given by $f := \frac{t}{\mathrm{vol}(\Sigma)} \cdot \mathbb{1}_\Sigma$.*

*Proof.* It is clear that $\hat{f} := \frac{t}{\mathrm{vol}(\Sigma)} \cdot \mathbb{1}_\Sigma$ satisfies the constraints of the problem (3.2) defining $I_\Omega^{\mathrm{TV}}(t)$ as soon as $t \le \mathrm{vol}(\Sigma)$, which ensures $0 \le \hat{f} \le \mathbb{1}_\Sigma \le \mathbb{1}_\Omega$. Hence, $I_\Omega^{\mathrm{TV}}(t) \le h_1(\Omega)t$. On the other hand, using the co-area formula (2.5), if $f$ is any competitor for the problem (3.2) then

$$\mathrm{TV}(f) = \int_0^{+\infty} \mathrm{area}(\partial\{f \ge s\})ds$$

$$= \int_0^{+\infty} \mathrm{vol}(\{f \ge s\}) \cdot \underbrace{\frac{\mathrm{area}(\partial\{f \ge s\})}{\mathrm{vol}(\{f \ge s\})}}_{\ge h_1(\Omega) \text{ by definition}} ds$$

$$\ge h_1(\Omega) \int_0^{+\infty} \mathrm{vol}(\{f \ge s\})ds = h_1(\Omega) \int_{\mathbb{R}^d} f(x)dx = h_1(\Omega)t.$$

Hence, for $t \le \mathrm{vol}(C)$, we have $I_\Omega^{\mathrm{TV}}(t) = h_1(\Omega)t$. ∎

If $\Omega$ is not convex, then there may exist more than one Cheeger set in $\Omega$ [34, §4]. In particular, in light of the proof of Proposition 3.6, there may be more than one minimizer of (3.2). Since we use an interior point solver in our experiments, the presence of multiple solutions for non-convex boundaries explains why our solutions appear fuzzy in this case.

Considering Proposition 3.6, a byproduct of our numerical method detailed below is the ability to compute Cheeger constants and Cheeger sets, at least when they are unique. This problem was tackled numerically in [13], where the authors propose a method to compute the largest Cheeger set of a given set $\Omega$. Compared to our numerical results, theirs are sharper; they also can characterize precisely which set is selected. On the other hand, we largely are interested in the optimal objective *value*, i.e., the Cheeger constant, which does not depend on the fuzziness of the minimizer $f$; moreover, the link to the Cheeger problem—and hence the numerical method of [13]—is lost as soon as $t$ becomes large. We leave as future work a detailed exploration of the relationship between these two techniques.

## 3.4. Structure of the minimizers.
As we have seen in the previous section, at least for small $t$ one can choose a solution $f$ of the problem (3.2) that is proportional to an indicator function, i.e., that takes only two values. Due to the potential non-uniqueness of Cheeger sets, however, it is possible to construct optimal functions $f$ that take on infinitely many values, even for small $t$. Up to the selection of appropriate optimal minimizers for the TV profile, however, we can prove the following:

**Proposition 3.7.** *There exists a family $(f_t)_{t \in [0,1]}$ such that:*
- *For any $t \in [0,1]$, the function $f_t \in L^1(\mathbb{R}^n)$ satisfies $0 \le f_t \le \mathbb{1}_\Omega$, $\int_{\mathbb{R}^n} f_t(x)\,dx = t$ and $\mathrm{TV}(f_t) = I_\Omega^{\mathrm{TV}}(t)$.*
- *For any $t \in [0,1]$, there exist $v_t \in (0,1)$ such that $f_t$ takes its values in $\{0, v_t, 1\}$.*
- *For a.e. $x \in \Omega$, the function $t \to f_t(x)$ is increasing.*

*Proof.* Let $S \subseteq [0, \mathrm{vol}(\Omega)]$ be the union of $\{0, \mathrm{vol}(\Omega)\}$ and the set of $t$ such that $I_\Omega^{\mathrm{TV}}$ is not an affine function in a neighborhood of $t$; the latter is closed in $[0, \mathrm{vol}(\Omega)]$. Since $I_\Omega$ is lower semi-continuous and $I_\Omega^{\mathrm{TV}}$ is its lower convex envelope, for any $t \in [0, \mathrm{vol}(\Omega)]$, if $I_\Omega^{\mathrm{TV}}(t) < I_\Omega(t)$ then $I_\Omega^{\mathrm{TV}}$ is affine in a neighborhood of $t$. In particular, $I_\Omega^{\mathrm{TV}}(t) = I_\Omega(t)$ for any point $t \in S$.

Let $t \le s$ with $t, s \in S$. By the previous remark, there exist subsets $\Sigma_t$ and $\Sigma_s$ of $\Omega$, with respective volumes $t$ and $s$, such that $\mathrm{area}(\partial \Sigma_t) = I_\Omega(t) = I_\Omega^{\mathrm{TV}}(t) = \mathrm{TV}(\mathbb{1}_{\Sigma_t})$; a similar identity holds for $\Sigma_s$. Let us consider $\Sigma_+ := \Sigma_t \cup \Sigma_s$ and $\Sigma_- := \Sigma_t \cap \Sigma_s$; we denote $r := \mathrm{vol}(\Sigma_-)$. By the inclusion–exclusion principle, $\mathrm{vol}(\Sigma_+) = t + s - r$. On the other hand, using [1, Proposition 3.38], the perimeter of $\Sigma_+$ satisfies

$$\mathrm{area}(\partial \Sigma_+) \le \mathrm{area}(\partial \Sigma_t) + \mathrm{area}(\partial \Sigma_s) - \mathrm{area}(\partial \Sigma_-).$$

Using $\mathbb{1}_{\Sigma_+}$ as a competitor for the problem (3.2) defining $I_\Omega^{\mathrm{TV}}(t + s - r)$, and given that $\mathrm{area}(\partial \Sigma_-) \ge I_\Omega^{\mathrm{TV}}(r)$, we deduce that

$$(3.7) \qquad I_\Omega^{\mathrm{TV}}(t + s - r) \le I_\Omega^{\mathrm{TV}}(t) + I_\Omega^{\mathrm{TV}}(s) - I_\Omega^{\mathrm{TV}}(r).$$

Convex non-negative functions are super-additive, but the equation above indicates that $I_\Omega^{\mathrm{TV}}$ exhibits sub-additive behavior. More precisely, using Lemma C.1 proved in the appendix, we see that $I_\Omega^{\mathrm{TV}}$ is affine on $[r, t+s-r]$ as soon as $r < t$. Since $t, s \in S$, the only way to avoid a contradiction is if $r = t$, which implies $\Sigma_t \subseteq \Sigma_s$.

As a consequence, using $t = s$, we see that for any $t \in S$ there exists a unique $\Sigma_t \subseteq \Omega$ such that $I_\Omega^{\mathrm{TV}}(t) = \mathrm{TV}(\mathbb{1}_{\Sigma_t})$. Using $t \le s$ we deduce that the map $t \in S \to \Sigma_t$ is increasing w.r.t. inclusion. In particular, for a.e. $x \in \mathbb{R}^n$, the map $t \in S \to \mathbb{1}_{\Sigma_t}(x)$ is increasing.[1] For $t \in S$, we set $f_t := \mathbb{1}_{\Sigma_t}$; such a choice satisfies all the requirements of Proposition 3.7.

Since the set $[0, \mathrm{vol}(\Omega)] \backslash S$ is open, it can be decomposed as a countable union of open intervals. Let $(t, s)$ be one of such interval, i.e., $t, s \in S$ and $S \cap (t, s) = \emptyset$. The definition of $S$, helped by a connectivity argument and the continuity of $I_\Omega^{\mathrm{TV}}$, shows that $I_\Omega^{\mathrm{TV}}$ is affine on the segment $[t, s]$. On the other hand, for any $r \in [t, s]$, the function

$$(3.8) \qquad f_r := \frac{s-r}{s-t}\mathbb{1}_{\Sigma_t} + \frac{r-t}{s-t}\mathbb{1}_{\Sigma_s},$$

i.e., the convex combination of $\mathbb{1}_{\Sigma_t}$ and $\mathbb{1}_{\Sigma_s}$ whose total mass is $r$, satisfies $\int_\Omega f_r(x)\,dx = r$, $0 \le f_r \le 1_\Omega$ and, by convexity of the TV norm,

$$\mathrm{TV}(f_r) \le \frac{s-r}{s-t}I_\Omega^{\mathrm{TV}}(t) + \frac{r-t}{s-t}I_\Omega^{\mathrm{TV}}(t).$$

The r.h.s. of this equation is precisely the affine function joining $(t, I_\Omega^{\mathrm{TV}}(t))$ to $(s, I_\Omega^{\mathrm{TV}}(s))$, which implies that the inequality is an equality. Hence, the function $f_r$ satisfies $0 \le f_r \le \mathbb{1}_\Omega$, $\int_{\mathbb{R}^n} f_r(x)\,dx = t$, and $\mathrm{TV}(f_r) = I_\Omega^{\mathrm{TV}}(r)$.

As $\Sigma_t \subset \Sigma_s$, the map $r \in [t, s] \to f_r(x)$ is increasing for a.e. $x \in \mathbb{R}^n$. Moreover, the function $f_r$ takes the value 0 on $\Omega \backslash \Sigma_s$, the value 1 on $\Sigma_t$, and only a third value, namely $(r-t)/(s-t)$ on the set $\Sigma_s \backslash \Sigma_t$. In conclusion, on every open interval within $[0, \mathrm{vol}(\Omega)] \backslash S$ we define $f_r$ by (3.8), and such a choice satisfies all the requirement of Proposition 3.7. ∎

### 3.5. Isoperimetric inequality.
So far, we have provided some propositions describing the shape of $I_\Omega^{\mathrm{TV}}$; our next task is to verify that it has properties in common with the isoperimetric profile that make it useful for evaluating compactness. Although our relaxation does not always recover the solution to the original problem, a key property is preserved:

**Proposition 3.8 (Isoperimetry).** *Suppose $B \subset \mathbb{R}^n$ is a ball whose volume matches $\mathrm{vol}(\Omega)$. Then, for all $t \in [0, \mathrm{vol}(\Omega)]$, we have $I_B^{\mathrm{TV}}(t) \le I_\Omega^{\mathrm{TV}}(t)$, and if the equality holds for some $t > 0$ then $\Omega$ is a ball.*

*Proof.* As recalled in the introduction, we already know that $I_B(t) \le I_\Omega(t)$ for any $t$. Hence, any convex function bounded by $I_B$ is also bounded by $I_\Omega$. Taking the supremum and recalling Proposition 3.5, we see that $I_B^{\mathrm{TV}}(t) \le I_\Omega^{\mathrm{TV}}(t)$.

Now assume that $I_B^{\mathrm{TV}}(t) = I_\Omega^{\mathrm{TV}}(t)$ for some $t \in (0, \mathrm{vol}(\Omega)]$. Taking into account the linearity of $I_B^{\mathrm{TV}}$ (Example 3.2) and the convexity of $I_\Omega^{\mathrm{TV}}$, we know that the slope of $I_B^{\mathrm{TV}}$ and $I_\Omega^{\mathrm{TV}}$ coincide at $t = 0$. Considering Proposition 3.6, this implies that the Cheeger constants of $B$ and $\Omega$ agree, which can only happen if $\Omega$ is a ball [35, Proposition 6.11]. ∎

---

[1] For this to hold we have to choose a precise representative of $\Sigma$. For instance, we choose the set of points in $\mathbb{R}^n$ such that the Lebesgue density of $\Sigma$ is equal to 1.

**4. Discretization and Optimization.** Having established theoretical properties of the total variation isoperimetric profile, our next step is to provide a discretization and algorithm for its approximation in practice.

**4.1. Discretization.** In our experiments, we assume that the shape $\Omega$ is expressed as an indicator function on a uniform grid; for example, when $n = 2$ we take as input an image $I$ with $I_{i,j} = 1$ inside $\Omega$ and $I_{i,j} = 0$ outside. We use the four-fold discretization $\nabla_p$ of the gradient operator $\nabla$ proposed in [15] to promote rotational invariance for our model. Namely, denote by $P$ the total number of grid points and by $\Delta x$ the grid size. Then, the linear operator $\nabla_p : \mathbb{R}^P \to \mathbb{R}^{P \times 4}$ is defined at any point $p = (i, j)$ as

$$(4.1) \qquad (\nabla_p f)_{i,j} = \frac{1}{\Delta x} \begin{pmatrix} f_{i+1,j} - f_{i,j} \\ f_{i,j+1} - f_{i,j} \\ f_{i+1,j+1} - f_{i,j+1} \\ f_{i+1,j+1} - f_{i+1,j} \end{pmatrix}.$$

Using this operator, the total variation of a function defined on the uniform grid is approximated by

$$\sum_p \|\nabla_p f\|_2 = \sum_{i,j} \|(\nabla_p f)_{i,j}\|_2$$

where $\| \cdot \|_2$ is the Euclidean norm in $\mathbb{R}^4$. After padding the grid with zeros on all four sides, we do not need to account for boundary conditions.

Our discretization of (3.2) becomes:

$$(4.2) \qquad I_\Omega^{\text{discretized}}(t) := \begin{cases} \inf_{f \in \mathbb{R}^P} & \sum_p \|\nabla_p f\|_2 \\ \text{subject to} & \mathbb{1}^\top f = t \\ & 0 \leq f \leq I \end{cases}$$

While careful proof of convergence in the limit of grid refinement is outside the scope of our discussion, we note that several analogous results exist for image processing models with similar structure [4, 31].

**4.2. Optimization.** Our problem (4.2) is a second-order cone program (SOCP) [9], for which there exist extremely efficient industrial solvers. We find that Mosek [2]—easily called using the CVX library [25, 26]—is effective for up to medium-scale instances and competitive with hand-designed algorithms. For scalability and simplicity, however, we can also derive a first-order algorithm based on the alternating direction method of multipliers (ADMM) [8], detailed below.

Suppose $z \in [0, 1]^{P'}$ is the restriction of the unknown variable $f$ in (4.2) to those pixels $p$ in the image $I$ that are nonzero; the remaining entries of $f$ must equal zero. We can write an equivalent formulation of problem (4.2) as

$$(4.3) \qquad I_\Omega^{\text{discretized}}(t) := \begin{cases} \inf_{z \in \mathbb{R}^{P'}} & \sum_p \|\mathcal{P}_p G z\|_2 \\ \text{subject to} & \mathbb{1}^\top z = t \\ & 0 \leq z \leq 1. \end{cases}$$

**Algorithm 4.1** ADMM TO SOLVE PROBLEM (4.4)

---

**function** TVPROFILEADMM
 $M \leftarrow \rho G^\top G + \tau \mathbb{1}\mathbb{1}^\top + \beta I_{P' \times P'}$      ▷ Can be factored once before iteration

 **for** $i = 1, 2, \ldots$ **do**
  *// Primal 1: z*
  $r \leftarrow G^\top(\rho x + y) + (t\tau - \lambda)\mathbb{1} + (\beta z' - q)$    ▷ Right-hand side of linear system
  $z \leftarrow M^{-1}r$        ▷ Can be accelerated using Cholesky factorization

  *// Primal 2: $(x, z')$*         ▷ Decouples over $x$ and $z'$
  *// Primal 2a: x*
  $H \leftarrow \text{RESHAPE}(Gz - y/\rho, 4 \times P)$     ▷ Isolate individual gradient vectors
  **for** $p = 1, 2, \ldots, P$ **do**        ▷ Per-column operation
   $c_p \leftarrow \text{MAX}(1 - 1/\rho \|H_p\|_2, 0)$      ▷ $H_p$ is column $p$ of $H$
   $X_p \leftarrow c_p H_p$
  $x \leftarrow \text{RESHAPE}(X, 4P \times 1)$      ▷ Unroll back to a vector

  *// Primal 2b: z'*
  $z' \leftarrow \text{CLAMP}(z + q/\beta, [0, 1])$

  *// Dual: $(y, \lambda, q)$*         ▷ Dual ascent
  $y \leftarrow y + \rho(x - Gz)$
  $\lambda \leftarrow \lambda + \tau(\mathbb{1}^\top z - t)$
  $q \leftarrow q + \beta(z - z')$

  *// Check for convergence and update $(\rho, \tau, \beta)$ here*   ▷ See [8, §3.3, §3.4] for details
 **return** $(x, z, z')$
**end function**

---

Here, $G \in \mathbb{R}^{4P \times P'}$ is the restriction of the gradient operator $\nabla$ in (4.1) to the nonzero pixels, and $\mathcal{P}_p \in \{0, 1\}^{4 \times 4P}$ extracts the four elements of the gradient relevant to pixel $p$.

To derive tractable ADMM iterations, we rewrite (4.3) in a somewhat counterintuitive form that leads to a tractable splitting:

$$(4.4) \qquad I_\Omega^{\text{discretized}}(t) := \begin{cases} \inf_{x, z, z'} & \sum_p \|\mathcal{P}_p x\|_2 + \chi(0 \le z' \le 1) \\ \text{subject to} & x - Gz = 0 & : y \\ & \mathbb{1}^\top z = t & : \lambda \\ & z - z' = 0. & : q \end{cases}$$

This expression copies $z$ into a second variable $z'$ and isolates the gradient $Gz$ as a third variable $x$; $\chi$ denotes a (convex) indicator function that equals $\infty$ any time the constraint is violated and 0 otherwise. The augmented Lagrangian of this optimization problem is

$$\Lambda(x, z, z'; y, \lambda, q) := \sum_p \|P_p x\|_2 + \chi(0 \le z' \le 1) + \frac{\rho}{2}\|x - Gz\|_2^2 + y^\top(x - Gz)$$
$$+ \frac{\tau}{2}(\mathbb{1}^\top z - t)^2 + \lambda(\mathbb{1}^\top z - t) + \frac{\beta}{2}\|z - z'\|_2^2 + q^\top(z - z').$$

Our two-block ADMM scheme cycles between three steps:

(Primal 1)     $z \leftarrow \min_z \Lambda(x, z, z'; y, \lambda, q)$     *Linear system*

(Primal 2)     $(x, z') \leftarrow \min_{(x, z')} \Lambda(x, z, z'; y, \lambda, q)$     *Closed-form*

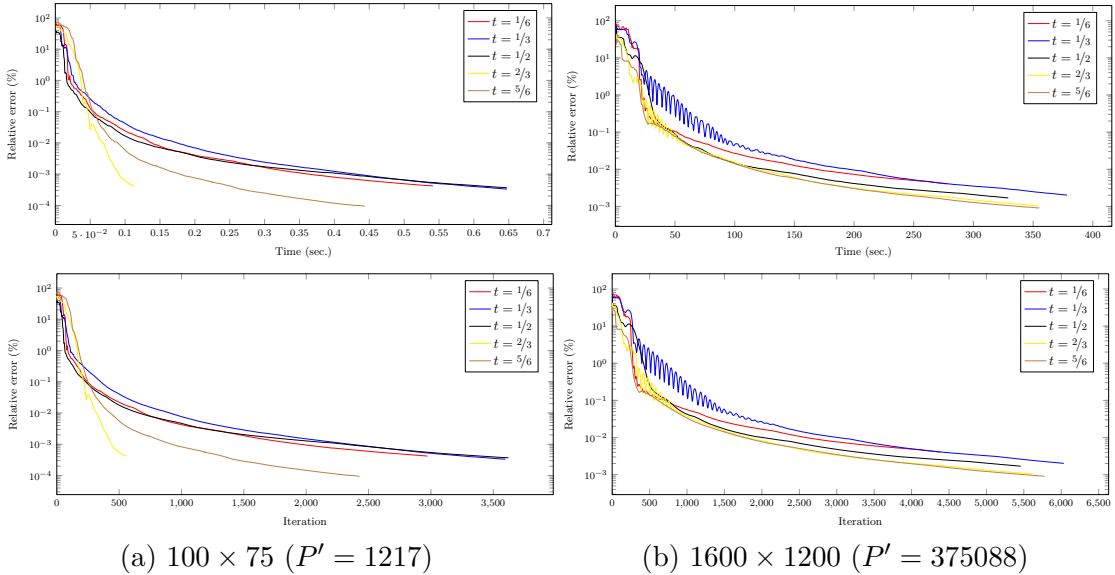

(a) $100 \times 75$ $(P' = 1217)$            (b) $1600 \times 1200$ $(P' = 375088)$

Figure 4: Error relative to ground truth for the ADMM algorithm presented in §4.2, as a function of time (top) and iteration (bottom). The image used for testing is district 12 in Figure 8. The different tests correspond to evenly-spaced $t$ values, $t \in \{^1/_6, \dots, ^5/_6\}$ for total area $A = 1$.

(Dual)     $\qquad\qquad\qquad y \leftarrow y + \rho(x - Gz) \qquad\qquad\qquad\qquad$ *Closed-form*

$\qquad\qquad\qquad\qquad\qquad \lambda \leftarrow \lambda + \tau(\mathbb{1}^\top z - t)$

$\qquad\qquad\qquad\qquad\qquad q \leftarrow q + \beta(z - z').$

Iteratively cycling these three steps is guaranteed to converge to the global optimum, as justified by the discussion in [8, §3.2.1]. Our splitting is carefully designed to make each iteration above computationally tractable; $z$ is obtained by solving a linear system of equations while the $(x, z')$ step decouples over these two variables and is solvable in closed-form (soft thresholding and clamping). The linear system for $z$ has the same matrix in each iteration, which can be prefactored. Algorithm 4.1 fills in the details of the steps.

We found this iterative technique to be the most efficient among the many possible methods for convex problems in the form (4.2); for instance, one alternative might be to use proximal splitting as described in [14, §3.2.4] and [16, §5]. It is worth noting, however, that we are in the worst-case situation for these algorithms: the terms in the most obvious proximal splitting are neither smooth nor strictly convex. Indeed, at least in the continuous case, the form of the optimal $f$ described in Proposition 3.7 indicates that typical solutions of our problem are at points of non-differentiability.

Figure 4 tests the efficiency of the ADMM algorithm in §4.2 on small-scale and large-scale examples. Our experiments involve computation of evenly-spaced samples from the TV profile of a nonconvex shape; the size of the image and number $P'$ of unknowns in the convex optimization problem are shown below the plots. Our experiments were carried out in Matlab 2018b, on an Ubuntu 16.04 machine with 64 GB memory and an Intel Xeon Gold 6136 CPU (3.00GHz).

For each example, we plot error of the ADMM variable $z$ relative to ground truth computed using Mosek [2] with relatively high precision ($\epsilon^{3/4}$, where $\epsilon \approx 2.22 \times 10^{-16}$ is machine precision); the relative error of the duplicated variable $z'$ converges similarly. While it is difficult to match Mosek's and ADMM's convergence criteria exactly, as a point of reference the ground truth computations in Mosek took an average of 0.726 seconds per small-scale example and 265.5 seconds per large-scale example to converge.

Our ADMM implementation is far from optimized, but it does include a few straightforward improvements to accelerate convergence and iteration time. In particular, we use the heuristic suggested in [8, §3.4.1] to adjust the ADMM penalty parameter up to 50 times during the optimization procedure; we check the heuristic once every 20 iterations. Sparse Cholesky factorization is used to prefactor the linear system for $z$ before iterations begin; this factorization has to be recomputed any time the penalty parameter changes.

Although ADMM has guaranteed linear convergence in theory, these plots verify that the convergence rates are practical even for large-scale examples, and that they are competitive with highly-optimized industrial solvers. Derivative discontinuities in the error plots largely correspond to automatic adjustments of the penalty parameter. Fairly high relative error ($\sim 0.1\%$) is tolerable for our target applications, which typically involve simply plotting the TV profile; from our experiments we can see that ADMM reaches this tolerance level quickly.

Source code including acceleration techniques above is provided in the GitHub repository accompanying this paper.[2]

## 5. Examples and Experiments.
We evaluate the value of total variation profiles through a number of experiments on synthetic shapes as well as American political districts. Our results not only confirm the theoretical properties explored in previous sections but also demonstrate that $I_\Omega^{\mathrm{TV}}(t)$ provides a fairly intuitive description of a shape at different length/area scales.

### 5.1. Synthetic Examples.
Figure 5 shows the optimized image $f$ for a variety of domains $\Omega$ and values of $t \in [0, 1]$. As suggested in §3.3, for small values of $t$ the images are simply rescalings of each other, concentrated in a compact subdomain of $\Omega$. For larger values of $t$, the optimization problem fills in progressively less compact regions within $\Omega$; at $t = 1.0$ the entire domain is filled.

The corresponding plots of $I_\Omega^{\mathrm{TV}}(t)$ are shown in Figure 6. In these and other plots, we scale the horizontal axis to the range $[0, 1]$ and the vertical axis by the perimeter of a circle whose area is vol($\Omega$). With this scaling, the profile of a circle is the diagonal line with slope 1. As would likely be expected, the text image is the least compact at all $t$ scales. Other more fine-grained information can be obtained by examining these plots, however. For instance, the square and hexagon are considered equally compact at nearly all $t \in [0, 1]$ except near $t \approx 1$. The spiral of circles connected by straight edges has similar compactness values for both edge thicknesses until $t \approx 0.9$, at which point the example with thinner edges distinguishes itself.

For values of $t$ less than 0.75, the dissected circle and the chain of circles connected by increasingly narrow bridges have almost identical total variation profiles. Examining Figure 5,

---

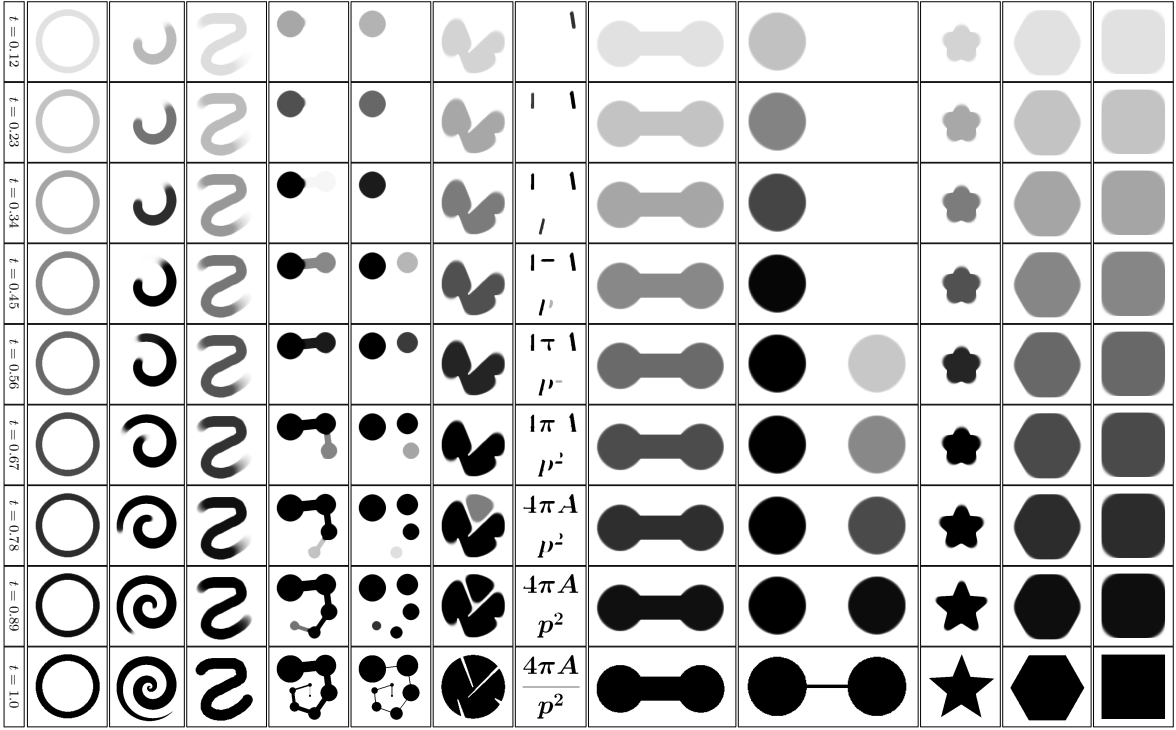

Figure 5: The optimal function $f$ as a function of $t$ for different shapes $\Omega \subseteq \mathbb{R}^2$; here, values of $f$ are scaled from zero (white) to one (black).

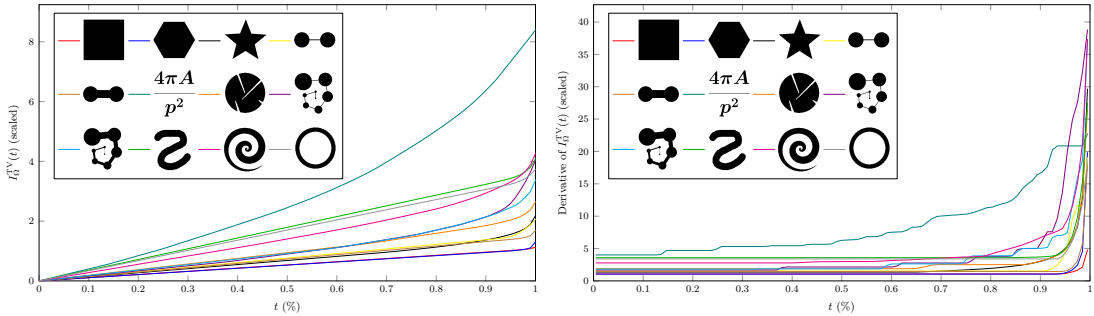

Figure 6: TV profiles (left) and their first derivatives (right) for the shapes in Figure 5.

we can see that at coarse scales both of the corresponding functions $f$ look like three relatively compact regions connected to each other. At finer scales, however, the dissected circle still appears as three compact connected regions whereas the three smallest circles in the chain are highly non-compact; this difference leads the two plots to diverge.

We can also observe that the donut shape and the S shape have similar total variation profiles at nearly all values of $t$, with the S being slightly less compact at all resolutions. Again examining the images in Figure 5, we can see that the S is highly symmetric and fills in almost uniformly; the donut is truly symmetric and fills in uniformly. Since the figures are of similar width and have similar perimeters, their total variation profiles are also similar.

Figure 6 also shows the first derivative of $I_\Omega^{\mathrm{TV}}(t)$, computed from the first plot using divided differences. As predicted by Proposition 3.6, all the plots have constant slope starting from $t = 0$. More surprisingly, the plots appear to have several flat regions, suggesting, together with Proposition 3.7 that the set $S$ that appears in the proof of this proposition is discrete. Note that a $t$ for which the derivative of $I_\Omega^{\mathrm{TV}}$ switches from one value to another corresponds to value for which the isoperimetric profile $I_\Omega(t)$ and the TV profile $I_\Omega^{\mathrm{TV}}(t)$ coincide.

### 5.2. Geographic Examples.

We next evaluate the total variation profile on a panel of real congressional districts defined by three different plans for North Carolina's districts for the U.S. House of Representatives. Following the 2010 Census, North Carolina was apportioned 13 seats and in 2011 the Republican legislature enacted the districts labelled **2011 districts**[3] which were in effect for the 2012 and 2014 congressional elections. After a long court process, this plan was ruled a racial gerrymander by the U.S. District Court for the Middle District of North Carolina in 2016 [27]; the U.S. Supreme Court affirmed the lower Court's decision in 2017 [19], deeming the 2011 plan unconstitutional.

Following these events, the Republican legislature enacted a new map, which took effect before the 2016 elections and is labelled **2016 districts**[4] in our figures. While more geographically compact, it was again challenged as a racial and partisan gerrymander and was struck down by the same U.S. District Court [17]. In 2019, the U.S. Supreme Court overruled the District Court, holding that partisan gerrymandering is not federally justiciable, and reinstated the 2016 plan [18].

Additionally, a panel of retired judges simulated a nonpartisan, independent redistricting commission and drew their own proposed plan, henceforth referred to as the **judges' plan**. An analysis performed by Common Cause North Carolina, using past election data from four election cycles, found that in contrast with the Legislature's enacted plans, which both have ten likely Republican and three likely Democratic districts, the judges' plan has six likely Republican, four likely Democratic, and three toss-up districts, which more closely matches North Carolina's statewide partisan vote shares [46].

We perform our analyses on North Carolina because all three plans are drawn with respect to the 2010 Census population data, reducing the number of confounding factors. Meaningful conclusions about relative compactness are difficult to draw when comparing plans from different states, since factors such as state boundaries, specific laws and rules governing the redistricting process, and number of seats can affect the compactness profile.

For this analysis, we extracted the polygonal boundary of each district, scaled to fit in a $250 \times 250$ bounding box; indicator functions of the polygon interiors sampled on this grid were used to approximate $I_\Omega^{\mathrm{TV}}(t)$. Plots of the TV profile for each of the 36 districts are shown in Figures 7(a–c). Examining the mean curves and standard deviation bands for the collection of districts in each plan in Figure 7(d), we can see that the 2011 plan appears less compact than the judges' plan, and the judges' plan appears less compact than the 2016 plan. The wide standard deviation band for the 2011 plan is strongly driven by the profile for the

---

twelfth district, which is extremely non-compact. This data reaffirms that non-compactness and gerrymandering are not equivalent: the 2016 map is more compact but less politically representative than the judges' plan.

The curves for the first and ninth districts in the 2011 plan cross around $t = 0.8$, with the ninth district appearing more compact at larger values of $t$. This occurs because the ninth district is consists of three small compact cores joined with thin connectors, while the first district has a large, fairly compact core on the northern side with several thin, snakelike pieces reaching in various directions. These pieces add a large amount of perimeter and a small amount of area, but the algorithm does not begin filling these appendages until roughly $t = 0.8$; at coarser resolutions, the region is fairly compact.

A similar observation can be made of the eighth district in the judges' plan, which consists of a compact core with a single arm stretching east. At coarse scales, the algorithm fills in this core; given that the other districts are also relatively compact, the curve for the eighth district sits in the middle of this range. Around $t = 0.75$, however, the curve bends upwards and crosses several others. Figure 10 reveals that the algorithm begins filling in this appendage around that time, whereas the algorithm is still filling the core in other districts.

Our analysis above shows the detailed information about compactness encoded in the TV profile. By examining and comparing TV profiles and the corresponding optimized functions $f$, we can *explain* quantitative scores for evaluating geometric quality.

**6. Extensions.** We briefly mention some potential extensions of $I_\Omega^{\mathrm{TV}}$ that may be of interest in different applications. In each case, we show that the extension is a small change of our basic convex optimization problem. Rather than adapting the technique in §4.2, for simplicity in this section we compute our examples using standard convex cone programming software [2]; our goal is to demonstrate qualitative aspects of these extensions empirically.

**6.1. Higher Dimensions.** The definition (3.2) and discretization (4.2) easily extend to dimensions $n > 2$. This extension allows us to evaluate the compactness of volumes embedded in $\mathbb{R}^3$ using a cubic lattice discretization of the domain.

Figure 11 shows example shapes computed using this volumetric version of the optimization problem for $I_\Omega^{\mathrm{TV}}(t)$. In this example, we represent $\Omega \subseteq \mathbb{R}^3$ as an indicator on a $100 \times 100 \times 100$ volume. While gathering enough samples to plot $I_\Omega^{\mathrm{TV}}(t)$ as a function of $t$ is prohibitively expensive, here we show the result of the optimization procedure at sparsely-sampled $t$ values; at each $t$ we render the level set of the unknown $f$ at the mean nonzero intensity. Once again, for small $t$ the shape remains fairly constant, as predicted by Proposition 3.6. Simple shapes like the torus remain preserved for most $t$'s, while more complex shapes like the humanoid start from a nearly-convex core when $t \ll 1$ and build up piece-by-piece.

**6.2. Accounting for Population Density.** Our profile (3.2) is defined to be purely geometric and is unaware of any measure on $\mathbb{R}^n$ other than the standard one. Given a distribution $\rho \in \mathrm{Prob}(\mathbb{R}^n)$—e.g. the population distribution of a state—we can extend our definition to evaluate compactness with respect to $\rho$:

$$(6.1) \qquad I_{\Omega,\rho}^{\mathrm{TV}}(t) := \left\{ \begin{array}{ll} \min_{f \in L^1(\mathbb{R}^n)} & \mathrm{TV}[f] \\ \text{subject to} & \int_{\mathbb{R}^n} f(x)\, d\rho(x) = t \\ & 0 \leq f \leq \mathbb{1}_\Omega. \end{array} \right.$$

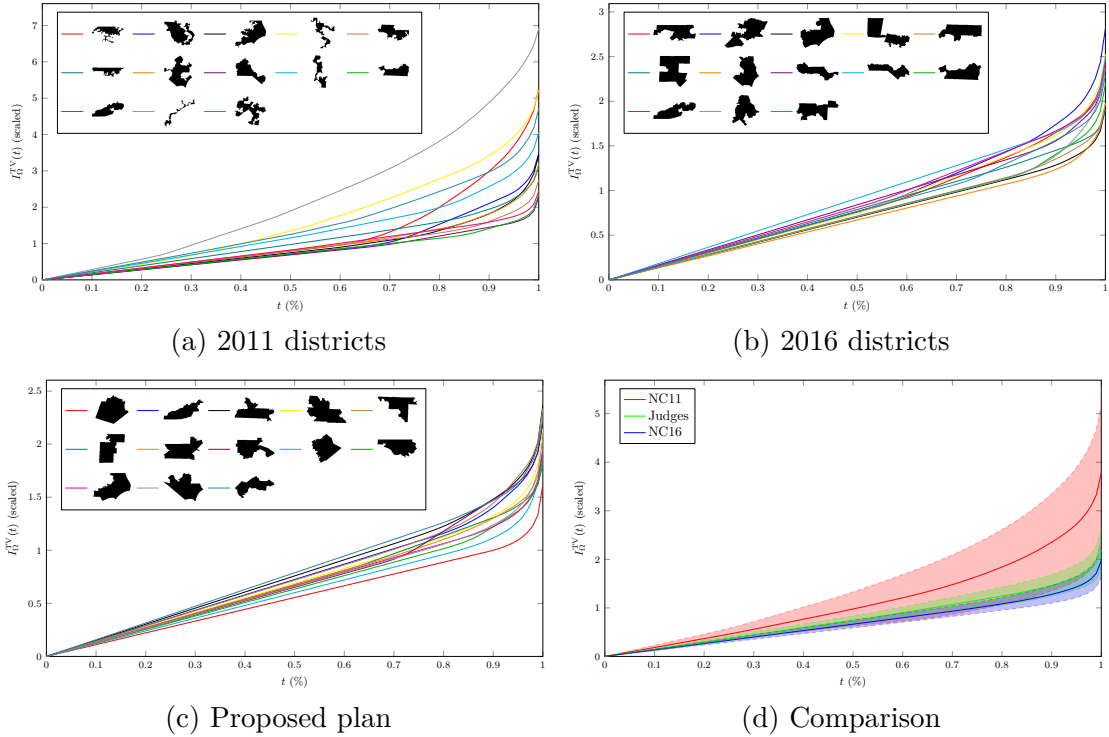

(a) 2011 districts

(b) 2016 districts

(c) Proposed plan

(d) Comparison

Figure 7: Total variation profiles for North Carolina districts enacted in 2011 (a; see Figure 8), enacted in 2016 (b; see Figure 9), and proposed by a nonpartisan panel of retired judges (c; see Figure 10). (d) The means and one-standard deviation bands for the total variation profiles for the three plans. In (a), (b), and (c) the districts are ordered sequentially from left to right and top to bottom to correspond with the numbering in Figures 8, 9, and 10.

A related formulation is investigated theoretically in [37], for the unrelaxed problem.

Figure 12 shows an experiment using this model. Here, we show the functions $f(x)$ computed on the same district shape, with three different density functions $\rho(x)$. Here we see how the choice of $\rho(x)$ can affect our assessment of compactness; it becomes less expensive to draw circles around densely-populated regions contained within the interior of the district, and the boundary of the district becomes less relevant because it is sparsely populated.

**6.3. Compactness on a Graph.** Even though the end result of political redistricting is a collection of geographic districts on a map, redistricting can often be described as a *graph theory* problem. In particular, districts commonly are built out of small geographic subunits, such as voting tabulation districts (VTDs) or Census blocks. Each VTD becomes a vertex in the graph, and two vertices are connected by an edge if and only if their corresponding VTDs are geographically adjacent.[5] We can then define a graph for each district as the one induced by the vertices corresponding to VTDs in that district. For this reason, it may be of interest

---

[5] *Adjacent* here is as according to *rook contiguity*, meaning that two VTDs must share a non-trivial segment of their borders to be considered adjacent. Two VTDs whose borders only share a single point are *not* considered adjacent in our construction.

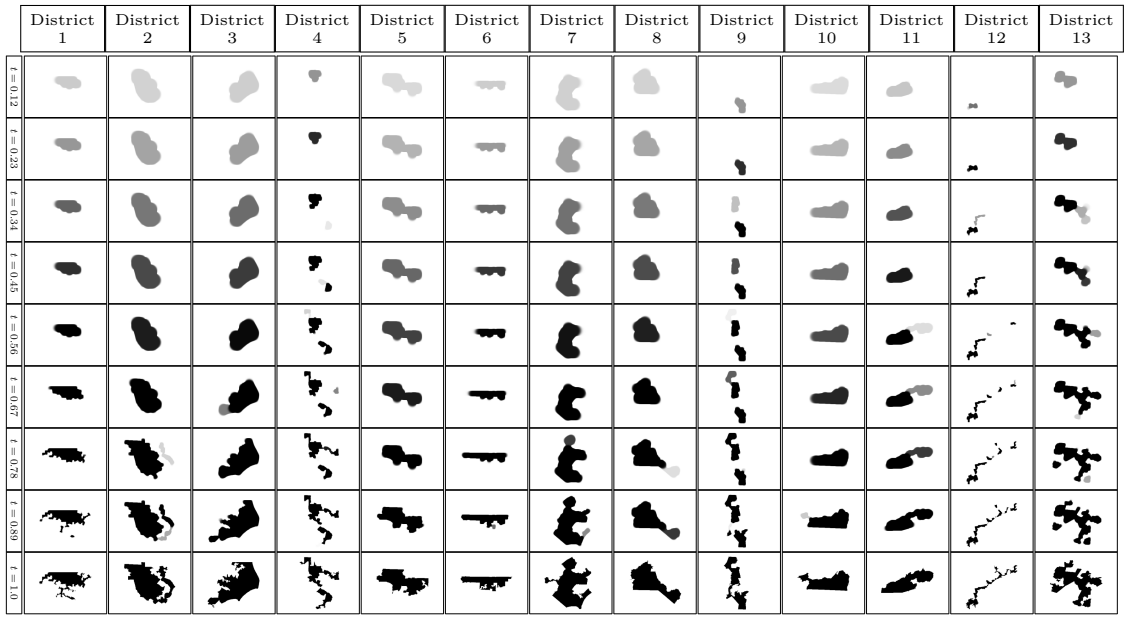

Figure 8: The optimal function $f$ as a function of $t$ for each of the 13 districts in the districting plan enacted in 2011; again, values of $f$ are scaled from zero (white) to one (black).

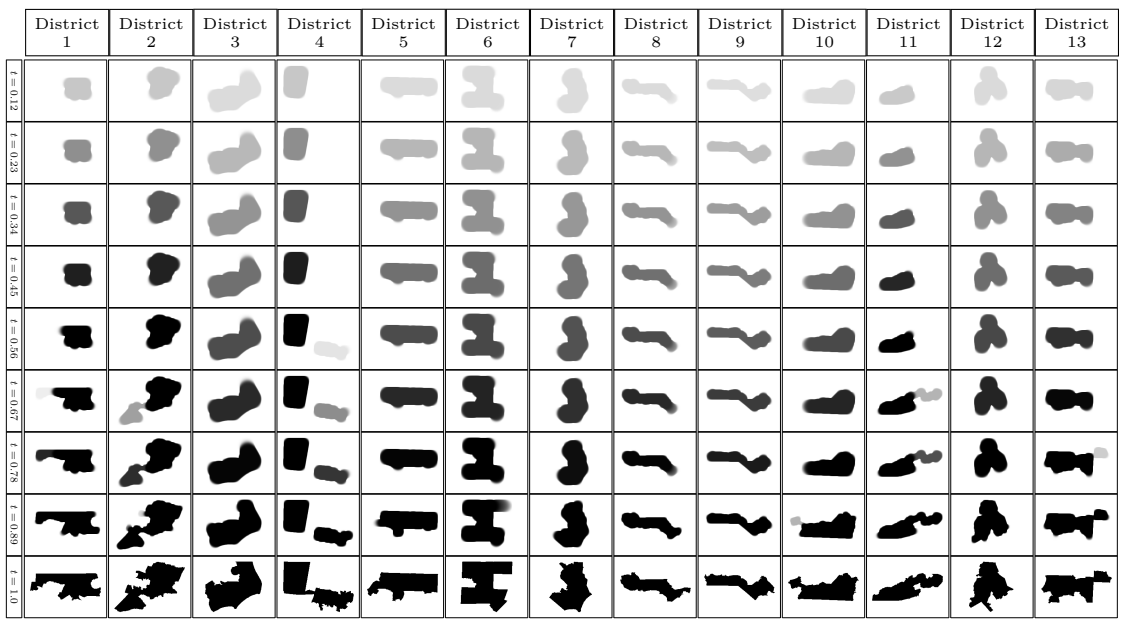

Figure 9: The optimal function $f$ as a function of $t$ for each of the 13 districts in the districting plan enacted in 2016; again, values of $f$ are scaled from zero (white) to one (black).

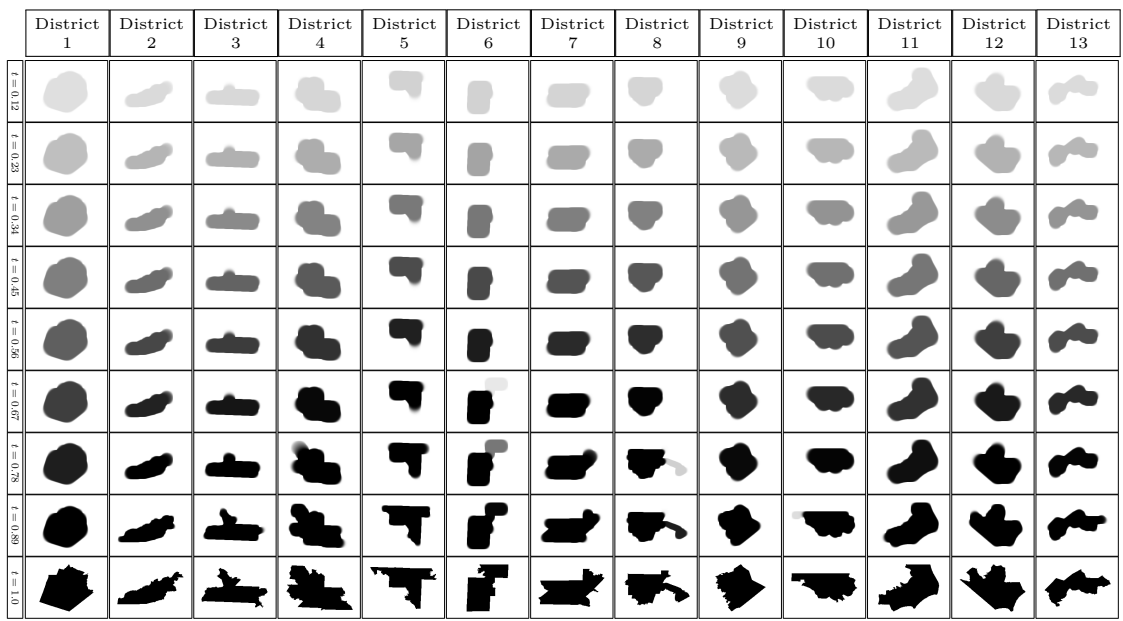

Figure 10: The optimal function $f$ as a function of $t$ for each of the 13 districts in the districting plan proposed by the nonpartisan panel of retired judges; again, values of $f$ are scaled from zero (white) to one (black).

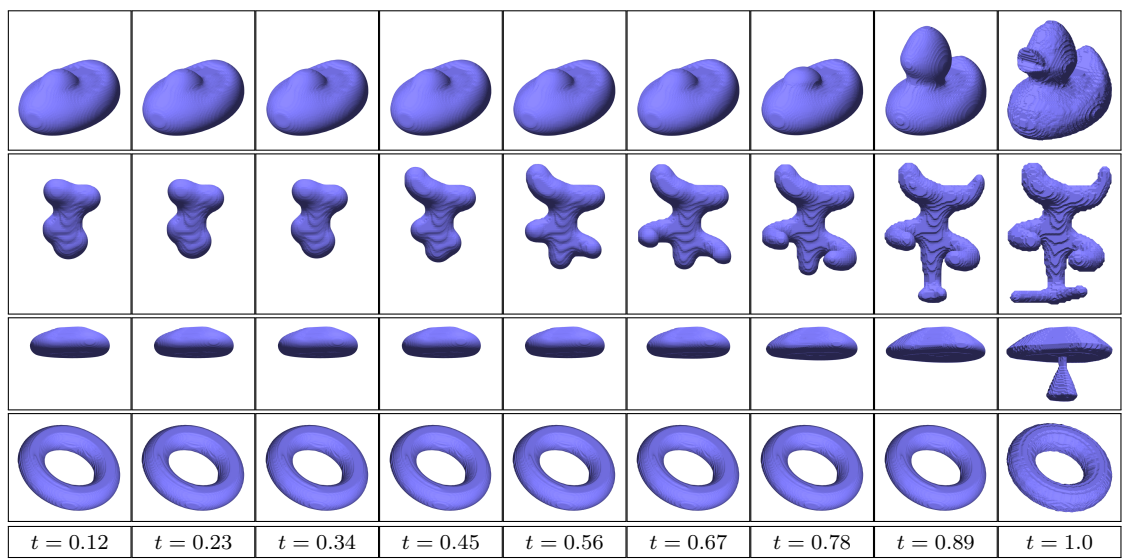

Figure 11: The identical optimization problem for TV profiles can be used for volumetric shapes in $\mathbb{R}^3$.

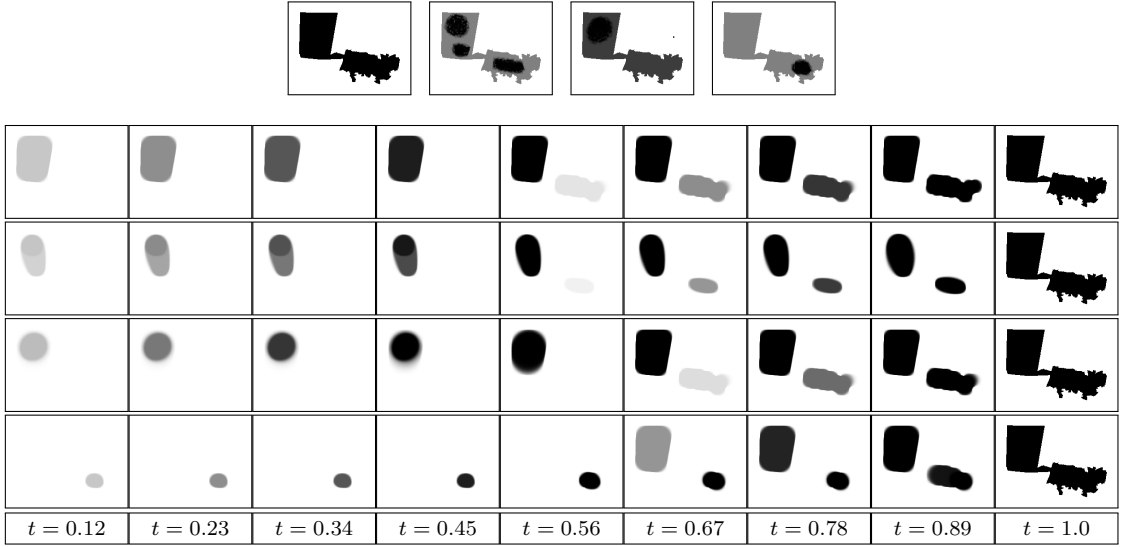

Figure 12: Fictional population distributions in the same district (top, left-to-right), and the frames of their corresponding weighted TV profiles (bottom, top-to-bottom).

to evaluate the compactness of a subset of graph vertices, bypassing embedding into $\mathbb{R}^n$.

Given a graph $G = (V, E)$, take $V_0 \subseteq V$ to be a subset of vertices representing a potential district. We can imitate (3.2) using graph-based constructions to define a TV profile of $V_0$:

$$(6.2) \qquad I_{V_0}^{\mathrm{TV}}(t) := \begin{cases} \min_{f \in \mathbb{R}^V} & \sum_{(v,w) \in E} |f(v) - f(w)| \\ \text{subject to} & \sum_{v \in V_0} f(v) = t|V_0| \\ & f(v) = 0 \ \forall v \notin V_0 \\ & f(v) \in [0, 1] \ \forall v \in V. \end{cases}$$

Note this definition is not the same as the version of the isoperimetric profile on a graph proposed in [38].

Figure 13 shows the graph-based total variation profile for a discretized version of North Carolina's 2011 districting plan, using the same underlying geography as the continuous experiments in §5.2. To perform the discretization, we use the U.S. Census Bureau's 2010 file of VTDs[6] and assign each VTD to the district in which it is contained.[7]

Figure 14 shows the corresponding plots of the TV profile for the discretized version of the 2011 North Carolina districting plan. As in the region-based profiles, the twelfth district is the least compact, and the fourth, tenth, and eleventh are all relatively compact compared to the rest of the districts. Beyond this connection, however, there are many points of difference between the two sets of profiles. One explanation is that areas with higher populations have more VTDs, and hence the graph-based model is in a sense aware of the population distribution. This makes it comparably less expensive to fill in densely-populated urban areas, which appear as compact cores under this measure. Conversely, geographically large rural zones have fewer VTDs and therefore appear less compact, as they are more expensive to fill.

---

[6] https://www.census.gov/geo/maps-data/data/tiger-line.html

[7] Since every VTD is not necessarily contained in a single congressional district, we use *areal interpolation* to assign VTDs to districts, where each is defined to belong to the district with which it shares the most area.

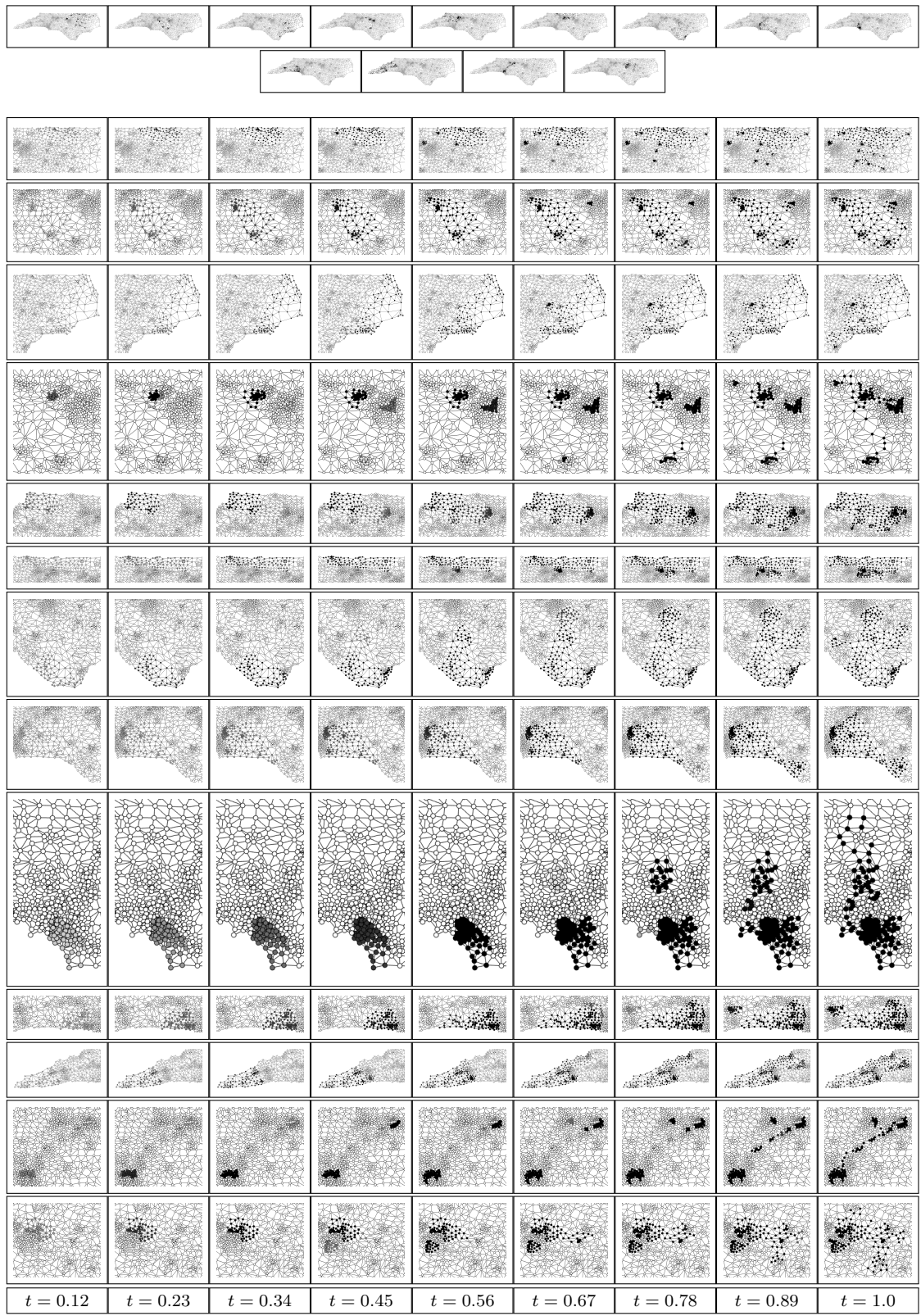

Figure 13: Total variation profile indicators for different subsets of vertices on a graph of North Carolina voting tabulation districts (VTDs). Images are cropped for each district; the top row shows how the districts are situated in the state. The algorithm is unaware of this embedding, and only sees an abstract representation of the graph.

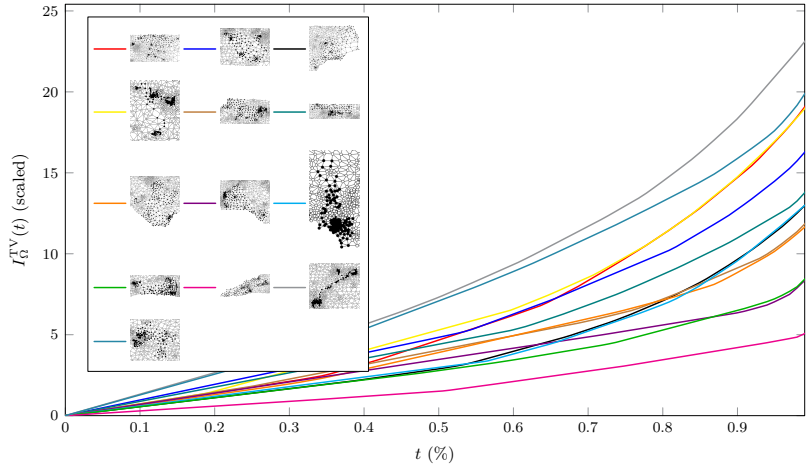

Figure 14: Graph-based TV profiles for the North Carolina districts in Figure 13.

To highlight this difference, recall that in the image-based model, most of the tenth district was filled-in fairly uniformly with the small chunk of Asheville in the northwest being filled at large values of $t$, but under the graph-based algorithm, the bulk of the district is filled in from east to west, beginning in the populous suburbs of Charlotte and spreading through the less heavily populated Western Foothills. Then, the portion of the district in Asheville is filled in before these two components are finally connected

## 7. Discussion and Conclusion.
Beyond their interest in the theory of geometry, isoperimetric profiles provide an intuitive, multiscale technique for evaluating the compactness of a shape $\Omega \subseteq \mathbb{R}^n$ that is not completely undermined by the instabilities of the isoperimetric quotient. By considering the entire plot, we gain a fairly complete description of a shape's level of contortion at different length/area scales. Furthermore, our TV-based relaxation of the classical profile admits simpler analysis using convex techniques.

Our initial work in proposing these profiles suggests several avenues for future mathematical, computational, and application-oriented research:

- A key consideration for our long-term target application of political redistricting is to evaluate methods for summarizing these plots and developing useful tests and benchmarks for existing and proposed districting plans. It will also be important to identify the most effective way to communicate the contents of the TV profile to an end user: non-mathematicians, including politicians and judges, must be able to interpret the results.
- Our formulation currently lacks a formal proof of convergence of our grid-based discretization in the limit of refinement, although this likely is a standard computation adjacent to results in total variation-based image denoising.
- Related to convergence, the graph-based model explored experimentally in §6.3 appears

to share many qualities with the TV profile for planar shapes. Exploiting the properties we have developed in the measure-theoretic case may lead to new insight into problems on graphs; ideally the two might be connected formally by examining convergence in the limit of refinement when the graphs are constructed by sampling planar regions.

- Inspired by the Lasserre hierarchy in semidefinite programming [32, 41], we might also ask whether there exist successively tighter convex approximations of the true isoperimetric profile, since our current relaxation is not tight (see Example 3.2); this question is intricately linked with Open Problem 2.1 above.
- Total variation is well-defined for functions on curved manifolds rather than the flat space $\mathbb{R}^n$; for instance, we could use our profile to evaluate the compactness of curved segments on a surface, e.g. in the case of geographic data accounting for the curvature and topography of the earth. To compute TV profiles in the presence of curvature, we easily could use a finite element formulation on triangulated surfaces, as suggested in [28].
- Our discussion focuses on total variation as an objective function, but we could attempt to generalize the construction of our profile by considering higher-order measurements popular in mathematical imaging like total generalized variation (TGV) [11].
- Finally, we could seek uniqueness results about our profile: Is it possible to encounter two shapes $\Omega$ with the same total variation profiles for all $t$?

These open questions aside, given our current analysis and optimization algorithm, the total variation profile is already a viable candidate for a nuanced and interpretable multiscale analysis of geometric compactness.

**Acknowledgments.** We thank Erik Demaine, Moon Duchin, Nestor Guillen, Simon Lam, Giorgio Saracco, Filippo Santambrogio, Lily Wang, and the participants of the 2018 Voting Rights Data Institute for their discussions and suggestions provided during the course of this research. We also thank Will Adler for providing areal interpolation code. J. Solomon acknowledges the generous support of the Prof. Amar G. Bose Research Grant, the MIT Research Support Committee ("Structured Optimization for Geometric Problems"), Army Research Office grant W911NF-12-R-0011 ("Smooth Modeling of Flows on Graphs"), National Science Foundation grant IIS-1838071 ("BIGDATA:F: Statistical and Computational Optimal Transport for Geometric Data Analysis"), Air Force Office of Scientific Research award FA9550-19-1-0319 ("Structured Assignment: Geometric Optimization Algorithms for Large-Scale Matching"), and an Amazon Research Award. Any opinions, findings, and conclusions or recommendations expressed in this material are those of the authors and do not necessarily reflect the views of these organizations.

**Appendix A. Regularity of the boundary $\partial\Omega$.** Throughout the text, the proofs about the properties of the TV profile work as long as

$$(A.1) \qquad\qquad I_\Omega^{\mathrm{TV}}(t) < +\infty \text{ for } t < \mathrm{vol}(\Omega).$$

Since $I_\Omega^{\mathrm{TV}}$ is convex (Proposition 3.4), this condition holds as soon as $I_\Omega^{\mathrm{TV}}(t)$ is finite for some $t$ arbitrarily close to $\mathrm{vol}(\Omega)$. We provide two cases where this can be shown to be true:

**Proposition A.1.** *Assume that* $\mathrm{area}(\partial\Omega) < +\infty$. *Then* (A.1) *holds.*

*Proof.* With this assumption, $I_\Omega^{\mathrm{TV}}(\mathrm{vol}(\Omega)) < +\infty$. Together with the convexity of $I_\Omega^{\mathrm{TV}}$ and $I_\Omega^{\mathrm{TV}}(0) = 0$, the conclusion clearly holds.                                    ∎

In the case where the boundary of $\Omega$ is fractal, we can still provide a criterion for (A.1).

**Proposition A.2.** *Assume that the topological boundary of $\Omega$ has zero (n-dimensional) Lebesgue measure. Then* (A.1) *holds.*

If the boundary of $\Omega$ is fractal, in the sense that its $d$-dimensional Hausdorff measure is finite for some $d \in (n-1, n)$, then we satisfy the assumptions of the proposition above.

*Proof.* Let $D$ be the distance function to $\mathbb{R}^n \backslash \Omega$. For $\delta < 0$, we set

$$\Sigma_\delta := \{x \in \Omega \text{ such that } D(x) \geq \delta\}.$$

Because of the assumption on the boundary of $\Omega$, $\mathrm{vol}(\Sigma_\delta) \to \mathrm{vol}(\Omega)$ as $\delta \to 0$. As the function $D$ is Lipschitz—and hence $\mathrm{TV}(D) < +\infty$—for a.e. $\delta$ we know that $\mathrm{area}(\partial \Sigma_\delta) < +\infty$. As a consequence, using $\mathbb{1}_{\Sigma_\delta}$ as a competitor, $I_\Omega^{\mathrm{TV}}(t) < +\infty$ for a.e. $t \in (0, \mathrm{vol}(\Omega))$. Together with the convexity of $I_\Omega^{\mathrm{TV}}$, we see that the conclusion holds.                              ∎

**Appendix B. Duality.**    Here we carefully justify the duality needed to derive (3.4). The proof relies on the Fenchel–Rockafellar duality theorem [12, Theorem 1.12]. For technical purposes (namely, compactness), we consider $B \subset \mathbb{R}^n$ a closed ball large enough to contain a neighborhood of $\Omega$. All the definitions considered in this paper do not change if we replace $\mathbb{R}^n$ by $B$.

Let $X = C(B)$ be defined as the space of continuous functions over $B$ valued in $\mathbb{R}$. Its dual $X^\star$ is the set of signed measures over $B$ [43, Theorem 2.14].[8] Define $F, G : X \to \mathbb{R} \cup \{+\infty\}$ as

$$F(\eta) := \int_\Omega \max(\eta(x), 0)\, dx,$$

and

$$G(\eta) := \begin{cases} -\lambda t & \text{if } \eta = \lambda - \nabla \cdot \phi \text{ with } \lambda \in \mathbb{R},\ \phi \in C_c^1(B \to \mathbb{R}^n) \text{ and } \|\phi\|_\infty \leq 1, \\ +\infty & \text{otherwise.} \end{cases}$$

$G$ is well-defined because $\lambda$ is uniquely defined from $\eta$ as the mean value of $\eta$, if the decomposition $\eta = \lambda - \nabla \cdot \phi$ holds. One can check that $F$ and $G$ are convex. Moreover, at the point $\eta = 0 \in X$, the function $F$ is continuous/finite and $G$ is finite.

Hence, we can apply the Fenchel–Rockafellar duality theorem [12, Theorem 1.12], which states that

$$(\mathrm{B.1}) \qquad \min_{f \in X} [F(x) + G(x)] = \sup_{\mu \in X^\star} [-F^\star(\mu) - G^\star(-\mu)],$$

where $F^\star, G^\star$ are the Legendre transforms of $F, G$. The l.h.s. can be written as

$$\begin{cases} \min_{\eta, \lambda, \phi} & -\lambda t + \int_\Omega \max(\eta(x), 0)\, dx \\ \text{subject to} & \|\phi\|_\infty \leq 1 \text{ and } \nabla \cdot \phi + \eta = \lambda, \end{cases}$$

---

[8]This would not hold if $B = \mathbb{R}^n$, justifying our introduction of $B$.

recovering the r.h.s. of (3.4) up to sign.

On the other hand, we can compute the Legendre transforms of $F$ and $G$. For $\mu \in X^\star$,

$$F^\star(\mu) = \sup_\eta \int_B \eta\mu - \int_\Omega \max(\eta(x), 0)\, dx.$$

If $\mu$ is negative somewhere, taking $\eta \leq 0$ at the same place, one gets larger and larger values for $F^\star(\mu)$ by integrating against $s\eta$, where $s \to +\infty$. On the other hand, if $\mu$ is larger than the Lebesgue measure restricted to $\Omega$ at some point, taking $\eta \geq 0$ and testing against $s\eta$, $s \to +\infty$, we reach the conclusion that $F^\star(\mu) = +\infty$. In particular, our argument implies that if $F^\star(\mu) < +\infty$ then $\mu$ has a positive density w.r.t. the restriction of the Lebesgue measure to $\Omega$ and that this density is between 0 and 1. In short,

$$F^\star(\mu) = \begin{cases} 0 & \text{if } d\mu(x) = f(x)\, dx \text{ with } 0 \leq f \leq \mathbb{1}_\Omega \text{ a.e.,} \\ +\infty & \text{otherwise.} \end{cases}$$

As far as $G$ is concerned,

$$\begin{aligned} G^\star(\mu) &= \sup_\eta \int_B \eta d\mu - G(\eta) \\ &= \sup\left\{ \int_B (\lambda - \nabla \cdot \phi) d\mu + \lambda t \; : \; \lambda \in \mathbb{R} \text{ and } \phi \in C_c^1(B \to \mathbb{R}^n), \|\phi\|_\infty \leq 1 \right\} \\ &= \sup_\phi \left\{ \int_B (\nabla \cdot \phi) d\mu \; : \; \phi \in C_c^1(B \to \mathbb{R}^n), \|\phi\|_\infty \leq 1 \right\} + \sup_{\lambda \in \mathbb{R}} \lambda \left( \int_B \mu + t \right). \end{aligned}$$

Provided that $\mu$ has a density $f$ w.r.t. the Lebesgue measure (and it is the case if $F^\star(\mu) < +\infty$), the first term is nothing else than $\mathrm{TV}(f)$, and the second is finite if and only if $\int_B f(x)dx = -t$. In short,

$$G^\star(f(x)dx) = \begin{cases} \mathrm{TV}(f) & \text{if } \int_B f(x)dx = -t, \\ +\infty & \text{otherwise.} \end{cases}$$

Hence, the r.h.s. of (B.1) reads exactly as $-I_\Omega^{\mathrm{TV}}(t)$, where the latter is defined in (3.2). In conclusion, (B.1) gives the equality between (3.2) and (3.4).

## Appendix C. (Generalized) sub-additivity for convex function.
It is known that convex functions have a tendency of being super-additive. Hence, if they also are sub-additive, they must behave like affine functions. This is the subject of the following lemma, which we used in the proof of Proposition 3.7:

**Lemma C.1.** *Let $F : [0, +\infty) \to [0, +\infty)$ be a nonnegative convex function and take $r < t \leq s$ such that*

$$(C.1) \qquad\qquad F(t + s - r) \leq F(t) + F(s) - F(r).$$

*Then the function $F$ is affine on $[r, t + s - r]$.*

*Proof.* Let $\Psi$ the function defined by $\Psi(x) := F(x) + F(s) - F(r) - F(x + s - r)$ (i.e. we consider $t$ as a variable). If $\Psi$ were differentiable we would easily conclude that it is non-increasing, but to prove it rigorously we use finite differences. If $r \leq x < y \leq s$,

$$(\text{C.2}) \qquad \frac{\Psi(y) - \Psi(x)}{y - x} = \frac{F(y) - F(x)}{y - x} - \frac{F(y + s - r) - F(x + s - r)}{y - x} \leq 0,$$

where the last inequality comes from the inequality $x \leq x + s - r$ and the convexity of $F$.

On the one hand, $\Psi(r) = 0$, and on the other hand $\Psi(t) \geq 0$ thanks to (C.1). As $\Psi$ is non-increasing, it is identically 0 on $[r, t]$. Plugging back this information into (C.2) and using the fact that the finite difference quotients are non-decreasing thanks to convexity of $F$, we deduce that

$$\frac{F(y) - F(x)}{y - x} = \frac{F(y + s - r) - F(x + s - r)}{y - x} = \text{const.}$$

for any $r \leq x < y \leq t$. Hence the function $F$ is affine on $[r, t]$ and $[s, t + s - r]$ with the same slope, which implies, by convexity of $F$, that $F$ is affine on $[r, t + s - r]$. ∎

Note that we needed $r < t$ in the proof; otherwise we could not plug such an $x < y$ in (C.2).